\documentclass[]{amsart}
\usepackage{amssymb}
\usepackage{amsmath}

\author{Raf Cluckers$^*$}
\address{Katholieke Universit Leuven,
Departement wiskunde, Celestijnenlaan 200B, B-3001 Leu\-ven,
Bel\-gium. Current address: \'Ecole Normale Sup\'erieure,
D\'epartement de ma\-th\'e\-ma\-ti\-ques et applications, 45 rue
d'Ulm, 75230 Paris Cedex 05, France}
\email{raf.cluckers@wis.kuleuven.ac.be}
\urladdr{http://www.wis.kuleuven.ac.be/algebra/Raf/}
\thanks{$^*$ The author is a postdoctoral  fellow of the Fund for Scientific
Research - Flanders (Belgium) (F.W.O.)}

\date{}

\theoremstyle{plain}
\newtheorem{theorem}{Theorem}[section]
\newtheorem{proposition}[theorem]{Proposition}

\newtheorem{lemma}[theorem]{Lemma}
\newtheorem{cor}[theorem]{Corollary}

\theoremstyle{definition}
\newtheorem{definition}[theorem]{Definition}

\newtheorem{remark}[theorem]{Remark}

\newcommand{\Z}{\mathbb{Z}}

\newcommand{\C}{\mathbb{C}}
\newcommand{\N}{\mathbb{N}}
\newcommand{\Q}{\mathbb{Q}}
\newcommand{\Reg}{\mathrm{Reg}}
\newcommand{\Supp}{\mathrm{Supp\, }}

\newcommand{\aaa}{\alpha}
\newcommand{\bbb}{\beta}

\newcommand{\ddd}{\delta}
\newcommand{\lP}{\lambda P_n}

\DeclareMathOperator{\sq}{\square}

\newcommand{\cC}{{\mathcal C}}
\newcommand{\cP}{\mathcal{P}}
\newcommand{\abs}[1]{\lvert#1\rvert}

\title[Decay rates of exponential sums]{Multivariate Igusa theory:\\
 Decay rates of exponential sums}
\begin{document}

\begin{abstract}
We obtain general estimates for exponential integrals of the form
\[
E_f(y)=\int_{\Z_{p}^{n}}\psi(\sum_{j=1}^r y_j f_j(x))|dx|,
\]
where the $f_j$ are restricted power series over $\Q_p$,
$y_j\in\Q_p$, and $\psi$ a nontrivial additive character on
$\Q_p$. We prove that if $(f_1,\ldots,f_r)$ is a dominant map,
then $|E_f(y)|<c|y|^{\alpha}$ for some $c>0$ and $\alpha<0$,
uniform in $y$, where $|y|=\max(|y_i|)_i$. In fact, we obtain
similar estimates for a much bigger class of exponential
integrals. To prove these estimates we introduce a new method to
study exponential sums, namely, we use the theory of $p$-adic
subanalytic sets and $p$-adic integration techniques based on
$p$-adic cell decomposition. We compare our results to some elementarily
 obtained explicit bounds for $E_f$ with $f_j$ polynomials. 
\end{abstract}

\maketitle

\noindent

\section{Introduction}
For $f=(f_1,\ldots,f_r)$ an $r$-tuple of restricted power series
over $\Q_p$ in the variables $x=(x_1,\ldots, x_n)$ and for $y\in
\Q_p^r$, we consider the exponential integral
\[
E_f(y)=\int_{\Z_{p}^{n}}\psi(y \cdot f(x) )|dx|,
\]
where $\psi$ is a nontrivial additive character on $\Q_p$, $|dx|$
denotes the normalized Haar measure on $\Q_p^n$, and $y \cdot
f(x)=\sum_jy_jf_j(x)$.
\par
With $|y|=\max(|y_i|)_i$ and $\ll$  the Vinogradov symbol, we
obtain the following general upper bounds:
\begin{theorem}\label{thm:decayexp:A}
If $f(\Z_p^n)$ has nonempty interior in $\Q_p^r$, there exists a
real number $\alpha<0$ such that
$$
E_f(y)\ll \min\{\abs{y}^{\alpha},1\}.
$$
\end{theorem}

In his book \cite{Igusa3} of 1978, J.~Igusa proves Theorem
\ref{thm:decayexp:A} in the case that $r=1$ with $f=f_1$ a
nonconstant homogeneous polynomial, and he formulates the problem
of generalizing this to the case of $r>1$. In this case Igusa is
also able to give an explicit $\alpha<0$ in terms of the numerical
data of an embedded resolution of $f$. By a very fine analysis of
embedded resolutions of $f$, Lichtin \cite{Lichtin2} is able to
prove Theorem \ref{thm:decayexp:A} in the case of a dominant map
of $r=2$ polynomials, where he also gives an explicit $\alpha<0$
in terms of the geometry of $f$. At present, these proofs seem to
be difficult to be generalized to the case of $r>2$ polynomials.
In the case of one nonconstant polynomial $f=f_1$, Theorem
\ref{thm:decayexp:A} can be proven by elementary methods, see for
instance the work of Chubarikov \cite{Chuba} and Loxton
\cite{Lox}. In the last section we show how the results of
\cite{Chuba} can be used to derive Theorem \ref{thm:decayexp:A},
when $f_1,\ldots,f_r$ are polynomials, in a very short way (even
with explicit upper bounds and weaker suppositions). We indicate
there why the situation for analytic maps is more difficult.

\subsection{}
In this paper we present a new technique to study exponential
integrals of a general nature, namely, by studying rather general
$p$-adic integrals by means of $p$-adic cell decomposition and the
theory of subanalytic sets. Examples of such general exponential
integrals are given below in this introduction. 
These techniques are also used in other contexts, for example, by
Denef \cite{Denef} to prove the rationality of the
Serre-Poincar\'e series associated to the $p$-adic points on a
variety.



\subsection{}
For readers not familiar with $p$-adic integration, we indicate
how $E_f(y)$ can be understood as an exponential sum. In the case
that the $f_i$ are restricted power series over $\Z_p$,
$\psi(x)=\exp(2\pi i (x\bmod\Z_p))$ (abbreviated by $\exp(2\pi i
x)$), and
 \begin{equation}\label{tupley}
 y=(\frac{u_1}{p^{m}},\ldots,\frac{u_r}{p^{m}}),
 \end{equation}
with $u_i$ integers satisfying $(u_1,\ldots,u_r,p)=1$, $m\geq0$,
we can write
$$
E_f(y)=\frac{1}{p^{mn}}\sum_{x\in (\Z_p/p^{m})^{n}}\exp(2\pi i
\frac{\sum_{j=1}^r u_j f_j(x)}{p^{m}}).
$$
Note that for general $y'\in \Q_p^r$ there can always be found a
tuple $y$ of the form (\ref{tupley}) such that $E_f(y')=E_f(y)$.
Theorem \ref{thm:decayexp:A} then says that $|E_f(y)|$ can be
bounded by $cp^{m\alpha}$ for some $c>0$ and $\alpha<0$, uniform
in $y$.

\subsection{}

We use the notion of subanalytic sets as in \cite{DvdD} and the
recent notion of subanalytic constructible functions as in
\cite{Ccell} (see below for the definitions).
\par
Let $G:\Q_p^r\to \Q$ be an integrable subanalytic constructible
function, and let $G^*(y):=\int_{\Q_p^r}G(x)\psi(x\cdot y)|dx|$ be
its Fourier transform. We obtain the following general upper
bounds:
\begin{theorem}\label{thm:decaysimple:A}
There exists a real number $\alpha<0$ such that
$
G^*(y) \ll \min\{\abs{y}^{\alpha},1\}.
$
\end{theorem}


\subsection{} We indicate how Theorem \ref{thm:decayexp:A} follows
from Theorem \ref{thm:decaysimple:A}. It is well known that,
whenever $f(\Z_p^n)$ has nonempty interior in $\Z_p^r$, $E_f$ is
the Fourier transform of an integrable function $F_f:\Q_p^r\to \Q$
(see \cite{Igusa3} or \cite{Weil}). We prove that we can take
$F_f$ to be a subanalytic constructible function (see Theorem
\ref{prop:fourier} below). Theorem \ref{thm:decayexp:A} then
follows immediately from Theorem \ref{thm:decaysimple:A}.
\par
In fact, a similar reasoning leads to the following much more
general result:
\begin{theorem}\label{thm:decayexp:B}
If $f:\Q_p^n\to \Q_p^r$ is a subanalytic map and
$\phi:\Q_p^n\to\Q$ is an integrable subanalytic constructible
function such that the support of $\phi$ is contained in
$f^{-1}(\mbox{Regular values of }f)\cup A$, with $A$ a set of
measure zero, then there exists a real number $\alpha<0$ such that
$E_{\phi,f}(y)\ll \min\{ |y|^{\alpha},1\}$, with
\[
E_{\phi,f}(y):=\int_{\Q_p^{n}}\phi(x)\psi(y\cdot f(x))|dx| .
\]
\end{theorem}

\par
We end section \ref{deca} with an open question about what happens
if $f$ is analytic but no longer subanalytic. All results of the
paper also hold for finite field extensions of $\Q_p$.

\par
Possible applications, or, possible subjects for future research,
lie in the search for candidate exponents $\alpha$ of Theorem
\ref{thm:decayexp:A} using the numerical data of a (parameterized)
resolution of singularities of the family $\sum_{i=1}^r u_if_i$
with parameters $u_i$ (if such resolution exists); as noted
before, candidate exponents can be found in this way when $r=1$,
see \cite{DenefBour}. Also, one can try to establish, under
similar conditions as in \cite{Igusa3}, an analytic analogue of
the Poisson summation formula considered by Igusa \cite{Igusa3}.

\subsection{Notation and terminology}
We fix a $p$-adic field  $K$ (i.e. $[K:\Q_p]$ is finite) and write
$R$ for the valuation ring of $K$, $\pi_0$ for a uniformizer of
$R$, and $q$ for the cardinality of the residue field. For $x\in
K$, $v(x)\in\Z\cup\{\infty\}$ denotes the $p$-adic valuation of
$x$ and $|x|=q^{-v(x)}$ the $p$-adic norm. We write $P_n$ for the
collection of $n$-th powers in $K^\times=K\setminus\{0\}$, $n>0$,
and $\lambda P_n=\{\lambda x\mid x\in P_n\}$ for $\lambda\in K$.
Let $\psi$ be a nontrivial additive character on $K$. We write $
x\cdot y=x_1y_1+\ldots+x_ny_n$ for $x,y\in K^n$, $n>0$.

The Vinogradov symbol $\ll$ has its usual meaning, namely that for
complex valued functions $f$ and $g$ with $g$ taking non-negative
real values $f\ll g$ means $|f|\leq c g $ for some constant $c$.

\par
A restricted analytic function $R^n\to K$ is an analytic function,
given by a single restricted power series over $K$ in $n$
variables (by definition, this is a power series over $K$ which
converges on $R^n$). We extend each restricted analytic function
$R^n\to K$ to a function $K^n\to K$ by putting it zero outside
$R^n$. A key notion is the following:
\begin{definition}
A subset of $K^n$ is called (globally) \emph{subanalytic} if it
can be obtained in finitely many steps by taking finite unions,
intersections, complements and linear projections of zero loci of
polynomials and of zero loci of restricted analytic functions in
$K^{n+e}$, $e\geq0$. A function $f:X\subset K^{m}\rightarrow
K^{n}$ is called subanalytic if its graph is a subanalytic set.
\end{definition}
We recall a basic result on subanalytic sets:
\begin{proposition}[\cite{DvdD}, Proposition (3.29)]\label{thm:basic:sub:an}
Let $X\subset K^n$  be a subanalytic set and $f:X\to K$ a
subanalytic function. Then there exists a finite partition of $X$
into $p$-adic submanifolds $A_j$ of $K^n$ such that the
restriction of $f$ to each $A_j$ is analytic and such that each
$A_j$ is subanalytic.
\end{proposition}
We refer to \cite{Ccell}, \cite{DenefBord, Denef1}, \cite{DvdD},
and \cite{vdDHM} for the theory of subanalytic sets.

\section{Cell decomposition and $p$-adic integration}\label{sec:decay:cell:decomp}
Cell decomposition is well suited to describe piecewise several
kinds of $p$-adic maps, for example,  polynomials maps, restricted
analytic maps, subanalytic constructible functions, and so on. It
allows one to partition the domain of such functions into $p$-adic
manifolds of a simple form, called cells, and to obtain on each of
these cells a nice description of the way the function depends on
a specific special variable (for an example of such an
application, see Lemma \ref{prop:descrip:simple}). By induction
one gets a nice description of the function with respect to the
other variables.
\par
Cells are defined by induction on the number of variables:
 \begin{definition}\label{def::cell}
A cell $A\subset K$ is a (nonempty) set of the form
 \begin{equation}
\{t\in K\mid |\aaa|\sq_1 |t-c|\sq_2 |\bbb|,\
  t-c\in \lP\},
\end{equation}
with constants $n>0$, $\lambda,c\in K$, $\aaa,\bbb\in K^\times$,
and $\square_i$ either $<$ or no condition. A cell $A\subset
K^{m+1}$, $m\geq0$, is a set  of the form
 \begin{equation}\label{Eq:cell:decay}
 \begin{array}{ll}
\{(x,t)\in K^{m+1}\mid
 &
 x\in D, \  |\aaa(x)|\sq_1 |t-c(
 x)|\sq_2 |\bbb(x)|,\\
 &
  t-c(x)\in \lP\},
  \end{array}
\end{equation}
 with $(x,t)=(x_1,\ldots,
x_m,t)$, $n>0$, $\lambda\in K$, $D=\pi_m(A)$ a cell where $\pi_m$
is the projection $K^{m+1}\to K^m$, subanalytic functions
$\aaa,\bbb:K^m\to K^\times$ and $c:K^m\to K$, and $\square_i$
either $<$ or no condition, such that the functions $\aaa,\bbb$,
and $c$ are analytic on $D$. We call $c$ the center of the cell
$A$ and $\lambda P_n$ the coset of $A$.
 \end{definition}
Note that a cell is either the graph of an analytic function
defined on $D$ (namely if $\lambda=0$), or, for each $x\in D$, the
fiber $A_x=\{t\mid (x,t)\in A\}$ is a nonempty open (if
$\lambda\not=0$).

\begin{theorem}[$p$-adic cell decomposition, \cite{Ccell}]\label{thm:CellDecomp}
Let $X\subset K^{m+1}$ be a subanalytic set and $f_j:X\to K$
subanalytic functions for $j=1,\ldots,r$. Then there exists a
finite partition of $X$ into cells $A_i$ with center $c_i$ and
coset $\lambda_i P_{n_i}$ such that
 \begin{equation}
 |f_j(x,t)|=
 |\ddd_{ij}(x)|\cdot|(t-c_i(x))^{a_{ij}}\lambda_i^{-a_{ij}}|^\frac{1}{n_i},\quad
 \mbox{ for each }(x,t)\in A_i,
 \end{equation}
with $(x,t)=(x_1,\ldots, x_m,t)$, integers $a_{ij}$, and
$\ddd_{ij}:K^m\to K$ subanalytic functions, analytic on
$\pi_m(A_i)$, $j=1,\ldots,r$. If $\lambda_i=0$, we use the
convention that $a_{ij}=0$.
 \end{theorem}
Theorem \ref{thm:CellDecomp} is a generalisation of cell
decomposition for polynomial maps by Denef \cite{Denef},
\cite{Denef2}. Recently, in \cite{Denef1} and \cite{Ccell}, cell
decomposition has been used to study parametrized integrals, as
follows.
\begin{definition}\label{basic algebra's}
For each subanalytic set $X$, we let $\cC(X)$ be the $\Q$-algebra
generated by the functions $|h|$ and $v(h)$ for all subanalytic
functions $h:X\to K^\times$. We call $G\in\cC(X)$ a
\emph{subanalytic constructible function} on $X$.
\par
To any function $G$ in $\cC(K^{m+n})$, $m,n\geq 0$, we associate a
function $I_m(G):K^m\to \Q$ by putting
 \begin{equation}\label{I_l}
I_m(G)(x)= \int\limits_{K^n}G(x,y)|dy|
 \end{equation}
if the function $y\mapsto G(x,y)$ is absolutely integrable for all
$x\in K^m$, and by putting $I_m(G)(x)=0$ otherwise.
\end{definition}
 \begin{theorem}[Basic Theorem on $p$-adic Analytic Integrals \cite{Ccell}]\label{thm:basic}
For any function $G\in\cC(K^{m+n})$, the function $I_m(G)$ is in
$\cC(K^{m})$.
 \end{theorem}

\begin{lemma}\label{prop:descrip:simple}
Let $X\subset K^{m+1}$ be a subanalytic set and let $G_j$ be
functions in $\cC(X)$ in the variables $(x_1,\ldots,x_m,t)$ for
$j=1,\ldots,r$. Then there exists a finite partition of $X$ into
cells $A_i$ with center $c_i$ and coset $\lambda_i P_{n_i}$ such
that each restriction $G_j|_{A_i}$ is a finite sum of functions of
the form
\[
|(t-c_i(x))^{a}\lambda^{-a}|^\frac{1}{n_i}v(t-c_i(x))^{s}h(x),
\]
where $h:K^m\to\Q$ is a subanalytic constructible function, and
$s\geq 0$ and $a$ are integers. Also, for any function $G\in
\cC(K^n)$ there exists a closed subanalytic set $A\subset K^n$ of
measure zero such that $G$ is locally constant on $K^n\setminus
A$.
\end{lemma}
\begin{proof}
The description is immediate from Theorem \ref{thm:CellDecomp} and
the definitions. The statement about $G\in\cC(K^n)$ follows from
Proposition \ref{thm:basic:sub:an} and the definitions.
\end{proof}
The following corollary is immediate.
\begin{cor}\label{lemma:decay} Let $G$ be in $\cC(K)$.
Suppose that if $|y|$ tends to $\infty$ then $G(y)$ converges to
zero. Then there exists a real number $\alpha<0$ such that
$G(y)\ll |y|^{\alpha}$. 
\end{cor}



We prove the following addendum to Theorem \ref{thm:basic}:
\begin{proposition}\label{prop:integrable}
Let $G$ in $\cC(K^{r+n})$ be such that $G(x,\cdot):K^n\to \Q$ is
integrable for almost all $x\in K^r$. Then, there exists a
function $F\in \cC(K^r)$ such that for all $x\in K^r\setminus B$,
with $B$ a subanalytic set of measure zero, one has
$$
F(x)= \int\limits_{K^n}G(x,y)|dy|.
$$
\end{proposition}
\begin{proof}
By induction and by Fubini's theorem it is enough to treat the
case $n=1$.

By Lemma \ref{prop:descrip:simple}, we can partition $K^{r+1}$
into cells $A$ with center $c$ and coset $\lambda P_m$ such that
$G{|_A}$ is a finite sum of functions of the form
\begin{equation}\label{eq:term:H:A}
 H(x,y)=|(y-c(
x))^a\lambda^{-a}|^\frac{1}{m}v(y-c(x))^{s}h(x),
\end{equation}
where $h:K^{r}\to \Q$ is a subanalytic constructible function,
and $s\geq0$ and $a$  are integers. \\

\textbf{Claim 1. }\textit{ Possibly after refining the partition,
we can assure that for each $A$ either the projection
$A':=\pi_{r}(A)\subset K^{r}$ has zero measure, or we can write
$G{|_A}$ as a sum of terms $H$ of the form (\ref{eq:term:H:A})
such that the function $H(x,\cdot)$ is integrable over
$A_{x}:=\{y\mid (x,y)\in A\}$ for all $x\in A'$.
}\\

First we prove the claim. By partitioning further, we may suppose
that either $v(y-c)$ is constant on $A$, or, it takes infinitely
many values on $A$, and in the case that $v(y-c)$ is constant on
$A$, we may assume that $a=s=0$. Regroup the terms with the same
exponents $(a,s)$, by summing up the respective functions $h$.

By the description (\ref{eq:term:H:A}) of $H$ and by the
definition of cells, the fact that the function
 \[
H(x,\cdot):A_{x}\to \Q: y\mapsto H(x,y)
 \]
is integrable over $A_{x}$ only depends on the exponents $(a,s)$,
on the fact whether $h(x)$ is zero or not, and on the particular
form of the cell $A_{x}$. Also, if  terms $H_1,\ldots,H_k$ have
different exponents $(a_i,s_i)$, then they have a different
asymptotical behavior for $y$ going to $c(x)$ with $x$ fixed, and
hence, if their sum is integrable over $A_{x}$, then each $H_i$ is
integrable over $A_{x}$.

Suppose now that $A$ has nonempty interior. 
Let $H$ be a term with exponents $(a,s)$ and function $h$ as in
(\ref{eq:term:H:A}). Then, either $h(x)$ is almost everywhere
zero, or, there exists by Lemma \ref{prop:descrip:simple} a
nonempty open $U\subset A'$ such that $h(x)$ is constant and
nonzero on $U$. If there exists such nonempty $U$, then, by the
above discussion, the term $H(x,\cdot)$ is integrable over $A_{x}$
for each $x\in U$ and hence for each $x\in A'$. If $h(x)$ is
almost everywhere zero, then we can, by partitioning $A'$ further
using Lemma \ref{prop:descrip:simple}, reduce to the case that
$A'$ has zero measure or $h(x)$ is identically zero on $A'$, in
which case we can skip the term $H$.
 This proves the claim.\\

Suppose that the statements of the claim are fulfilled for our
partition of $K^{r+1}$ into cells. Let $\cP$ be the set of cells
$A$ such that $\pi_{r}(A)$ has measure zero. Put
$B:=\cup_{A\in\cP}\pi_{r}(A)$ and $C:=\cup_{A\in\cP} A$. Let $G'$
be the constructible function $G(1-\chi_C)$ where $\chi_C$ is the
characteristic function of $C$. Then, $B$ has measure zero in
$K^r$ and $G'$ satisfies
$$
\int\limits_{K}G(x,y)|dy|=\int\limits_{K}G'(x,y)|dy|
$$
for all $x\in K^r\setminus B$ and $G'(x,\cdot)$ is integrable for
all $x\in K^r$. Putting $F:=I_r(G')$, an application of Theorem
\ref{thm:basic} ends the proof.
\end{proof}

\section{Exponential sums as Fourier tranforms}\label{sec:local}
We fix a nontrivial additive character $\psi$ on $K$. For
$\phi\in\cC(K^n)$ an integrable function, for $f:K^n\to K^r$ a
subanalytic function, and for $y\in K^r$, we consider the
exponential integral
 \[
E_{\phi,f}(y)=\int_{K^{n}}\phi(x)\psi(y\cdot f(x))|dx|.
 \]
We call $z\in K^r$ a \emph{regular value} of $f$ if $f^{-1}(z)$ is
nonempty, if $f$ is $C^1$ on a neighborhood of $f^{-1}(z)$, and if
the rank of the Jacobian matrix of $f$ is maximal at each point
$x\in f^{-1}(z)$. We denote the set of regular values of $f$ by
$\Reg_f$ and the support of $\phi$ by $\Supp \phi$ .

\begin{theorem}\label{prop:fourier}
Let $f:K^n\to K^r$ be a subanalytic function and let
$\phi\in\cC(K^n)$ be an integrable function satisfying $\Supp
\phi\subset f^{-1}(\Reg_f)\cup A$ with $A$ a set of measure zero.
Then there exists an integrable function $F_{\phi,f}$ in
$\cC(K^r)$ such that for any bounded continuous function
$G:K^r\to\C$ one has
$$
\int_{K^r}F_{\phi,f}(z)G(z)|dz|=\int_{K^n}\phi(x)G(f(x))|dx|,
$$
and hence, the following Fourier transformation formula holds:
 \[
E_{\phi,f}(y)=\int_{z\in K^{r}}F_{\phi,f}(z)\psi(z\cdot y)|dz|.
 \]
\end{theorem}

Theorem \ref{prop:fourier} is a generalisation of Corollary~1.8.2
in \cite{Denef1} by Denef which treats the case that the $f_i$ are
polynomials and $\phi$ is a Schwartz-Bruhat function. Igusa has
given an analogon of Theorem \ref{prop:fourier} in the case of
$r=1$ polynomial (cf. the asymptotic expansions of \cite{Igusa3}),
and Lichtin \cite{Lichtin} in the case of $r=2$ polynomials, both
in the case that $\phi$ is a Schwartz-Bruhat function. Igusa and
Lichtin also relate the asymptotic expansions to the numerical
data of an embedded resolution of $f$, the counterpart (however
not easily computable) of which would be here to apply cell
decomposition to get explicit asymptotic expansions for given $f$
and $\phi$.

 \par
Note that $F_{\phi,f}$ is determined up to a set of measure zero
by the universal property stated in the Theorem. The function
$F_{\phi_{\rm triv},f}$, with $\phi_{\rm triv}$ the characteristic
function of $R^n$ and $f$ a dominant polynomial mapping, is called
the local singular series of $f$ and plays an important role in
number theory, for example in the circle method.

\begin{proof}[Proof of Theorem \ref{prop:fourier}.]
Clearly $f^{-1}(\Reg_f)$ is subanalytic.
Without loss of generality we may assume that for all $x\in
f^{-1}(\Reg_f)$ one has
$$
J(x):=\det\left(\frac{\partial f_i}{\partial
x_j}(x)\right)_{i,j=1,\ldots,r}\not=0.
$$

By the inverse function theorem, Proposition
\ref{thm:basic:sub:an}, Theorem~(3.2) of \cite{DvdD} on the
existence of bounds, and the subanalytic selection Theorem (3.6)
of \cite{DvdD}, we may also suppose that
$$
T:f^{-1}(\Reg_f)\to K^n:x\mapsto y  = (f(x),x_{r+1},\ldots,x_n)
$$
is injective and a $C^1$ bijection onto its image with $C^1$
inverse. Applying the change of variables formula, we obtain
\begin{eqnarray*}
\int_{K^n}\phi(x)G(f(x))|dx| & = &  \int_{f^{-1}(\Reg_f)}\phi(x)G(f(x))|dx| \\
 & = &
 \int_{T(f^{-1}(\Reg_f))}\phi\circ T^{-1}(y)G(y_1,\ldots,y_r) | J\circ T^{-1}(y) |^{-1}
 |dy|.
\end{eqnarray*}
By Fubini's theorem and Proposition \ref{prop:integrable}, there
exists a function $F_{\varphi,f}$ in $\cC(K^r)$ with the property
that
$$
F_{\varphi,f}(y_1,\ldots,y_r)= \int_{K^{n-r}} \phi\circ
T^{-1}(y)\, | J\circ T^{-1}(y) |^{-1}
 |dy_{r+1}\wedge\ldots\wedge dy_n|,
$$
for almost all $(y_1,\ldots,y_r)\in K^r$, where we have extended
the integrand by zero to $K^{n-r}$. This function clearly
satisfies the requirements of the theorem.
\end{proof}

\section{Estimates for Fourier transforms}
For an integrable  function $G$ in $\cC(K^r)$ we write $G^*$ for
its Fourier transform
$$G^*:K^r\to\C:y\mapsto\int_{K^{r}}G(x)\psi(x\cdot y)|dx|.$$

The following is a generalisation of Theorem
\ref{thm:decaysimple:A}. 

\begin{theorem}\label{thm:decaysimple}
For each integrable $G\in\cC(K^r)$ there exists a real number
$\alpha<0$ such that $ G^*(y) \ll \min\{\abs{y}^{\alpha},1\}$.
\end{theorem}

\begin{proof}
For simplicity we suppose that $\psi(R)=1$ and $\psi(x)\not=1$ for
$x\not \in R$ (any other additive character is of the form
$x\mapsto\psi(ax)$ with $a\in K$). It is clear that $G^*(y)\ll 1$
since
\[
|G^*(y)|\leq \int_{K^{r}}|G(x)||dx|<\infty.
\]
Hence, it is enough to prove for $i=1,\ldots,r$ that
 $$G^*(y) \ll \abs{y_i}^{\alpha_i}
 $$
for some $\alpha_i<0$. We prove that $G^*(y) \ll
\abs{y_r}^{\alpha_r}$ for some $\alpha_r<0$. Write $x=(\hat
x,x_r)$ with $\hat
x=(x_1,\ldots,x_{r-1})$. 
By Lemma \ref{prop:descrip:simple}, we can partition $K^r$ into
cells $A$ with center $c$ and coset $\lambda P_m$ such that $G|_A$
is a finite sum of functions of the form
\begin{equation}\label{eq:term:H}
 H(x)=|(x_r-c(\hat
x))^a\lambda^{-a}|^\frac{1}{m}v(x_r-c(\hat x))^{s}h(\hat x),
\end{equation}
where $h:K^{r-1}\to \Q$ is a subanalytic constructible function,
and $s\geq0$ and $a$  are integers. \\

\textbf{Claim 2.}\textit{ Possibly after refining the partition,
we can assure that for each $A$ either the projection
$A':=\pi_{r-1}(A)\subset K^{r}$ has zero measure, or we can write
$G{|_A}$ as a sum of terms $H$ of the form (\ref{eq:term:H}) such
that $H$ is integrable over $A$ and $H(\hat x,\cdot)$ is
integrable over $A_{\hat x}:=\{x_r\mid (\hat x,x_r)\in A\}$ for
all $\hat x\in A'$. Moreover, doing so we can assure that each
such term $H$ does not change its sign on $A$. }\\

As this claim and its proof are similar to Claim 1 we will give
only an indication of its proof.

Partitioning further, we may suppose that $v(x_r-c(\hat x))$ does
not change its sign on $A$, and that it either takes only one
value on $A$ or infinitely many values. If $v(x_r-c(\hat x))$ only
takes one value on $A$, we may suppose that the exponents $a$ and
$s$ as in (\ref{eq:term:H}) are zero. Now apply Lemma
\ref{prop:descrip:simple} to each $h$ and to the norms of all the
subanalytic functions appearing in the description of the cells
$A$ in a similar way (in particular, make similar assumptions as
above). Do this inductively for each variable. This way, the claim
is reduced to a summation problem over (Presburger set of)
integers, which is easily solved (cf. the proof of Claim 1). This
proves the claim.


Fix a cell $A$ and a term $H$ as in the claim. The cell $A$ has by
definition the following form
\[
\begin{array}{ll}
A=\{x\mid & \hat x\in A',\ v(\alpha(\hat x))\sq_1 v(x_r-c(\hat
x))\sq_2 v(\beta(\hat x)),
\\
 &
 x_r-c(\hat x)\in\lambda P_m \},
\end{array}
\]
where $A'=\pi_{r-1}(A)$ is a cell, $\sq_i$ is $<$ or no condition,
and $\alpha,\beta:K^{r-1}\to K^\times$ and $c:K^{r-1}\to K$ are
subanalytic functions. We focus on a cell $A$ with nonempty
interior, in particular,  $\lambda\not=0$ and $A'$ has nonempty
interior. For $\hat x\in A'$ and $y\in K^r$, we denote by $I(\hat
x,y)$ the value
 \[
I(\hat x,y)=\int_{x_r\in A_{\hat x}} H(x)\, \psi( x\cdot y
)\,|dx_r|.
 \]
Let $\chi_{\lambda P_{m}}:K\to\Q$ be the characteristic function
of $\lambda P_{m}$ and write $\hat y=(y_{1},\ldots,y_{r-1})$. We
easily find that $I(\hat x,y)$ equals
\begin{equation}\label{eq:sum0}
\begin{array}{c}
\psi(\hat x\cdot \hat y+cy_r)\ h(\hat
x)|\lambda|^{-a/m}\sum\limits_{(\ref{summation})} q^{-ja/m}\,
j^{s} \int\limits_{v(x_r-c)=j} \chi_{\lambda
P_m}(x_r-c)\psi((x_r-c) y_{r})\,|dx_r|,
\end{array}
\end{equation}
where $c=c(\hat x)$ and the summation is over
 \begin{equation}\label{summation}
 \{j \mid v(\alpha(\hat x)) \sq_1 j \sq_2
v(\beta(\hat x))\}.
\end{equation}
\par
By Hensel's Lemma, there exists an integer $e$ such that all units
$u$ with $u\equiv 1\bmod \pi_0^e$ are $m$-th powers (here, $\pi_0$
is such that $v(\pi_0)=1$). Hence,
\[
\int_{v(u)=j}\chi_{\lambda P_m}(u)\psi(u y_{r})\,|du|
\]
is zero whenever $j+v(y_r)+e<0$ (since in this case one
essentially sums a nontrivial character over a finite group). By
consequence, the only terms contributing to the sum
(\ref{eq:sum0}) are those for which $-v(y_{r})-e\leq j$.

\par
We thus have
\begin{eqnarray}
|\int_{x\in A}H(x)\psi(x\cdot y)|dx|\, |
 & = &
 |\int_{\hat x\in A'} I(\hat x,y)|d\hat x|\,|\label{eq:int:I}
 \\
& \leq &
 \int_{B_{y_r}}|H(x)||dx|\label{eq:int:bound}
\end{eqnarray}
with 
 $ 
B_{y_r}= \{x\in K^r \mid x\in A,\ -v(y_r)-e \leq v(x_r-c(\hat x))
 \}.
 $ 

The integrability of $H$ over $A$, the fact that $H$ does not
change its sign on $A$, and Theorem \ref{thm:basic} imply that the
integral (\ref{eq:int:bound}), considered as a function in the
variable $y_r$, is in $\cC(K)$.
 \par
Next we prove that (\ref{eq:int:bound}) goes to zero when $|y_r|$
goes to infinity. First suppose that $A$ is contained in a compact
set. Since $B_{y_r}\subset A$, 
the measure of $B_{y_r}$, and hence also (\ref{eq:int:bound}),
goes to zero when $|y_r|$ tends to infinity. In the case that $A$
is not contained in a compact set, let $A_b$ be the intersection
of $A$ with $(\pi_0^bR)^r$, for $b<0$. Clearly each $A_b$ is
contained in a compact set. Also, for each $\varepsilon>0$, there
exists a $b_0$ such that for each $b<b_0$ and for each $y_r$ one
has $\int_{B_{y_r}\setminus A_b}|H(x)||dx|<\varepsilon$, by the
integrability of $H$ over $A$. By the previous discussion,
$\int_{B_{y_r}\cap A_b}|H(x)||dx|$, and hence also
(\ref{eq:int:bound}), goes to zero when $|y_r|$ goes to $\infty$.
\par
An application of Corollary \ref{lemma:decay} now finishes the
proof.

\end{proof}
\begin{remark}
The fact that $|G^*|$ in Theorem~\ref{thm:decaysimple} goes to
zero when $|y|$ goes to infinity also follows directly from the
Lemma of Riemann-Lebesgue in general Fourier analysis, 
cf.~the section on Fourier transforms in \cite{Weilgroup}.
However, to know this is not enough to apply Corollary
\ref{lemma:decay} as is done to finish the proof of
Theorem~\ref{thm:decaysimple} since in general $|G^*|$ is not
subanalytic constructible.
\end{remark}

\section{Decay rates of exponential sums}\label{deca} 
We use the notation of section \ref{sec:local} for $E_{\phi,f}$.
Combining Theorem \ref{thm:decaysimple} with the Fourier
transformation formula of Theorem \ref{prop:fourier}, we obtain
the following generalization of Theorem \ref{thm:decayexp:B}:
\begin{theorem}\label{thm:decayexp}
If $f:K^n\to K^r$ is a subanalytic map and $\phi\in\cC(K^n)$ is
integrable and satisfies $\Supp \phi\subset f^{-1}(\Reg_f)\cup A$
with $A$ a set of measure zero, then there exists a real number
$\alpha<0$ such that
\[
E_{\phi,f}(y) \ll \min\{ |y|^{\alpha},1\}.
\]
\end{theorem}
Combining this Theorem with the fact that the set of singular
points of a dominant polynomial mapping $K^n\to K^r$ (or a
dominant restricted analytic mapping $R^n\to K^r$) has measure
zero, we find:
\begin{cor}\label{cor:decayexp:B} If $f:K^n\to K^r$ is a dominant polynomial
mapping and if $\phi\in\cC(K^n)$ is integrable, then there exists
$\alpha<0$ such that
\[
E_{\phi,f}(y) \ll \min\{ |y|^{\alpha},1\}.
\]
The same conclusion holds for $E_{\phi,f}$ with $f:R^n\to K^r$ a
restricted analytic map, extended by zero to a map $K^n\to K^r$,
such that $f(R^n)$ has nonempty interior in $K^r$.
\end{cor}

\subsection{}
We end this section with an open question. Let
$f=(f_1,\ldots,f_r):K^n\to K^r$ be an analytic map given by $r$
power series $f_1,\ldots,f_r\in K[[x]]$ which converge on the
whole of $K^n$. Suppose that $\phi\in\cC(K^n)$ is integrable and
that $f(K^n)$ contains a nonempty open. The question is whether
there exists an $\alpha<0$ such that
$$
\int_{K^{n}}\phi(x)\psi(y\cdot f(x))|dx|\ll \min\{
|y|^{\alpha},1\}.
$$

\section{Polynomial mappings}
In this section we use elementary methods to deduce explicit upper
bounds for polynomial exponential sums. Theorem
\ref{prop:polyn:upper} below is of a different nature than our
main Theorem \ref{thm:decayexp} (and its proof is much more easy),
in the sense that it uses the degree of the polynomial mapping as
exponent in the upper bound. Such bound based on the degree would
give a trivial bound when naively adapted to the analytic case.
Similar problems occur when the explicit bounds of Loxton
\cite{Lox} are naively adapted to the analytic case. Since we use
a result of \cite{Chuba} formulated there for polynomials over
$\Z$, we will work over $\Q_p$.


For $g$ a polynomial in $\Q_p[x]$ with $x=(x_1,\ldots,x_n)$ let
$d_{j}(g)$ be the degree of $g$ with respect to the variable $x_j$
for $j=1,\ldots,n$, and let $e(g)$ be the minimum of the $p$-adic
orders of the coefficients of $g(x)-g(0)$. For
$f=(f_1,\ldots,f_r)$ a tuple of polynomials in $\Q_p[x]$ let
$d(f)$ be $\max_{ij}(d_{j}(f_i))$.

A function $\phi:\Q_p^n\to \Q$ is a Schwartz-Bruhat function if it
is locally constant and has compact support. In this section we
consider
\[
E_{\phi,f}(y)=\int_{\Q_{p}^{n}}\phi(x)\psi(y \cdot f(x))|dx|,
\]
with $f=(f_1,\ldots,f_r)$ a tuple of polynomials in $\Q_p[x]$,
$\phi:\Q_p^n\to \Q$ a Schwartz-Bruhat function, $\psi$ a
nontrivial additive character on $\Q_p$, and $y\in \Q_p^r$.

By elementary methods we easily deduce the following from work of
Chubarikov \cite{Chuba}.

\begin{theorem}\label{prop:polyn:upper}
Suppose that $f_1,\ldots, f_r$ are polynomials in
$x=(x_1,\ldots,x_n)$ over $\Q_p$ which satisfy that $\sum_i a_if_i
+a_0=0$  implies $a_i=0$ for $a_i\in \Q_p$ and $i=0,\ldots,n$. Let
$\phi:\Q_p^n\to \Q$ be a Schwartz-Bruhat function. Then, for any
$\varepsilon>0$, one has
$$
E_{\phi,f}(y)\ll \min\{\abs{y}^{\varepsilon-1/d(f)},1\}.
$$
 Moreover, for $y$ with $v(y)<0$, one has
$$
E_{\phi,f}(y)\ll (-v(y))^{n-1}\abs{y}^{-1/d(f)}.
$$
\end{theorem}

\begin{proof}
For simplicity we may assume that $\psi(\Z_p)=1$ and
$\psi(x)\not=1$ for $x\not \in \Z_p$ and that at least one
coefficient of $f_1(x)-f_1(0)$ has $p$-adic order $0$. Since
$\phi$ is a finite linear combination of characteristic functions
of compact balls, we may moreover assume that $\phi$ is $\phi_{\rm
triv}$, that is, the characteristic function of $\Z_p^n$.
Chubarikov \cite{Chuba}, Lemma 3, proves that for any polynomial
$g\in\Z[x]$ with $e(g)=0$, $d(g)\leq d$ for some $d\in\N$, and
each $z\in\Q_p$ with $v(z)<0$ one has
$$
|E_{\phi_{\rm triv},g}(z)| < c(d,n)(-v(z))^{n-1}|z|^{-1/d}
$$
with $c(d,n)$ a constant only depending on $d$ and $n$. 

Rewrite $E_{\phi,f}(y)$ as
\[
E'(z,u_1,\ldots,u_r)=\int_{\Z_{p}^{n}}\psi(z(u\cdot f(x)))|dx|,
\]
with $z\in \Q_p$, $u\in\Z_p^r$ with $|u|=1$, and
$y=(zu_1,\ldots,zu_r)$. For any such $u$, the number $d(u\cdot f)$
cannot exceed $d(f)$. By the compactness and completeness of
$\{u\in \Z_p^r\mid |u|=1\}$, also the number $e(u\cdot f)$ is
bounded uniformly in $u$, say, by $N$, since otherwise $\sum_i
a_if_i +a_0=0$ for some nontrivial $a_i\in \Q_p$.

One easily deduces from the mentioned result of \cite{Chuba} that
for $v(z)<-N$
 $$
  |E'(z,u_1,\ldots,u_r)|
< c(d(f),n) p^{N/d(f)} (-v(z))^{n-1} |z|^{-1/d(f)},
$$
with $c(d(f),n)$ as above. 
The Theorem follows.

\end{proof}
\begin{remark}
Note that from the proof of Theorem \ref{prop:polyn:upper} and
Lemma 3 of \cite{Chuba}, one can construct a (non optimal)
constant $c$, depending only on $\psi,\phi$, and $f$, such that
for each $y$ with $v(y)<0$
$$
|E_{\phi,f}(y)| < c(-v(y))^{n-1}\abs{y}^{-1/d(f)}.
$$
We leave the determination of the optimal $c$ for the future.
\end{remark}

\subsection*{Acknowledgment} I would like to thank J.~Denef, F.~Loeser, and B.~Lichtin for
interesting conversations on this and related subjects.

\bibliographystyle{amsplain}
\bibliography{anbib}

\providecommand{\bysame}{\leavevmode\hbox to3em{\hrulefill}\thinspace}
\providecommand{\MR}{\relax\ifhmode\unskip\space\fi MR }
\providecommand{\MRhref}[2]{%
  \href{http://www.ams.org/mathscinet-getitem?mr=#1}{#2}
}
\providecommand{\href}[2]{#2}
\begin{thebibliography}{10}

\bibitem{Chuba}
V.N. Chubarikov, \emph{Multiple rational trigonometric sums and multiple
  integrals}, Mat. Zametki \textbf{20} (1976), no.~1, 61--68, English transl.:
  Math Notes 20 (1976).

\bibitem{Ccell}
R.~Cluckers, \emph{Analytic $p$-adic cell decomposition and integrals}, Trans.
  Amer. Math. Soc. \textbf{356} (2004), no.~4, 1489 -- 1499.

\bibitem{Denef}
J.~Denef, \emph{The rationality of the {P}oincar\'e series associated to the
  $p$-adic points on a variety}, Inventiones Mathematicae \textbf{77} (1984),
  1--23.

\bibitem{Denef2}
\bysame, \emph{$p$-adic semialgebraic sets and cell decomposition}, Journal
  f{\"u}r die reine und angewandte Mathematik \textbf{369} (1986), 154--166.

\bibitem{DenefBord}
\bysame, \emph{Multiplicity of the poles of the poincar\'e series of a $p$-adic
  subanalytic set}, S\'em. Th. Nombres Bordeaux \textbf{43} (1987-1988), 1--8.

\bibitem{DenefBour}
\bysame, \emph{Report on {I}gusa's local zeta function}, S\'eminaire Bourbaki
  \textbf{Vol. 1990/91, Exp. No.730-744} (1991), 359--386, Ast\'erisque
  201-203.

\bibitem{Denef1}
\bysame, \emph{Arithmetic and geometric applications of quantifier elimination
  for valued fields}, MSRI Publications, vol.~39, pp.~173--198, Cambridge
  University Press, 2000.

\bibitem{DvdD}
J.~Denef and {L. van den} Dries, \emph{$p$-adic and real subanalytic sets},
  Annals of Mathematics \textbf{128} (1988), no.~1, 79--138.

\bibitem{vdDHM}
{L. van den} Dries, D.~Haskell, and D.~Macpherson, \emph{One-dimensional
  $p$-adic subanalytic sets}, Journal of the London Mathematical Society
  \textbf{59} (1999), no.~1, 1--20.

\bibitem{Igusa3}
J.~Igusa, \emph{Lectures on forms of higher degree (notes by {S}. {R}aghavan)},
  Lectures on mathematics and physics, Tata institute of fundamental research,
  vol.~59, Springer-Verlag, 1978.

\bibitem{Lichtin2}
B.~Lichtin, \emph{On a question of {I}gusa: towards a theory of several
  variable asymptotic expansions {II}}, preprint.

\bibitem{Lichtin}
\bysame, \emph{On a question of {I}gusa: towards a theory of several variable
  asymptotic expansions {I}}, Compositio Mathematica \textbf{120} (2000),
  no.~1, 25--82.

\bibitem{Lox}
J.~H. Loxton, \emph{Estimates for complete multiple exponential sums}, Acta
  Arith. \textbf{92} (2000), no.~3, 277--290.

\bibitem{Weilgroup}
A.~Weil, \emph{L'int\'egration dans les groupes topologiques et ses
  applications}, Publ. {I}nst. {M}ath. {C}lermont-{F}errand, Hermann \& Co.,
  Paris, 1940.

\bibitem{Weil}
\bysame, \emph{Sur la formule de {S}iegel dans la th{\'e}orie des groupes
  classiques}, Acta Mathematica \textbf{113} (1965), 1--87.

\end{thebibliography}

\end{document}